\documentclass{article}
\def\n{\noindent}
\begin{document}

\title{On the existence of universal series by trigonometric system }
\author{S.A.Episkoposian}

\date{e-mail: sergoep@ysu.am}

 \maketitle {
In this paper we prove the following: let $\omega(t)$ be a
continuous function, increasing in $[0,\infty)$ and
$\omega(+0)=0$. Then there exists a series of the form
 $$\sum_{k=-\infty}^\infty C_ke^{ikx} \ \  with \ \ \sum_{k=-\infty}^\infty C^2_k \omega(|C_k|)<\infty ,\ \ C_{-k}=\overline{C}_k, \eqno$$
with the following property: for each $\varepsilon>0$  a weighted
function $\mu(x), 0<\mu(x) \le1, \left | \{ x\in[0,2\pi]:
\mu(x)\not =1 \} \right | <\varepsilon $ can be constructed, so
that the series is universal in the weighted space
$L_\mu^1[0,2\pi]$ with respect to rearrangements.

} \vskip 4mm

{ \bf \S1. INTRODUCTION}

In 1932 F. Riesz(see [1], p. 655) proved that there exists a
function $f_0(x)\in L^1[0,2\pi]$ so that its Fourier series with
respect to the trigonometric system does not converge in $
L^1[0,2\pi]$. Consequently, there exist functions in the space
$L^1[0,2\pi]$ that cannot be represented by trigonometric series
in the metric of $L^1$.

Let $\mu(x)$ be a measurable on $[0,2\pi]$ function with $
0<\mu(x) \le1$, $ x\in[0,2\pi]$ and let $L_\mu^1[0,2\pi]$ be a
space of measurable functions $f(x),\ \ x\in [0,2\pi]$ with
$$\int_0^{2\pi} |f(x)| \mu(x) dx<\infty.$$

 K.Kazarian and R.Zink in [2] proved that there exist a weighted space $L_\mu^1[0,2\pi]$, such that for every function $f(x)$ in the space $L_\mu^1[0,2\pi]$ one can find a trigonometric series
 $\sum_{k=-\infty}^\infty C_ke^{ikx}$ that converges to $f(x)$ in the metric of $L_\mu^1[0,2\pi]$.

Moreover using other construction of weighted function $\mu(x)$
M.Grigorian proved the following result [3]

----------------------------------------------------------------

The author was supported in part by Grant- 08-83 from the Government of Armenia

{ {\bf AMS Classification} 2000 Primary 42A20 .}

\par\par\bigskip
\par\par\bigskip
\par\par\bigskip
\break

{\bf Theorem 1.} There exists a trigonometric series of the form
 $$\sum_{k=-\infty}^\infty{C_ke^{ikx} \ \  with \ \ \sum_{k=-\infty}^\infty \left | {C_k} \right|^q <\infty},\ \ \forall q>2\ \  C_{-k}=\overline
{C}_k\eqno(1.1)$$ with the following property: such that for any
number $\epsilon>0$  a weighted function $\mu(x)$, $ 0<\mu(x)
\le1, \left | \{ x\in[0,2\pi]: \mu(x)\not =1 \} \right | <\epsilon
$ can be constructed, so that the series (1.1) is universal in
$L_\mu^1[0,2\pi]$ with respect to subseries (see Definition 3).

Now we present the definitions of universal series:

{\bf Definition 1.}  A functional series
$$\sum_{k=1}^\infty f_k(x),\ \ f_k(x) \in L_\mu^1[0,2\pi] \eqno(1.2)$$
is said to be universal in weighted spaces $L_\mu^1[0,2\pi]$ with
respect to rearrangements, if for any function
 $f(x) \in L_\mu^1[0,2\pi]$ the members of (1.2) can be rearranged so that the obtained series $\displaystyle \sum_{k=1}^\infty f_{\sigma(k)}(x)$
converges to the function $f(x)$ in the metric  $L_\mu^1[0,2\pi]$,
i.e.
$$\lim_{n\to \infty} \int_0^{2\pi} \left| \sum _{k=1}^n f_{ \sigma (k)}(x)-f(x) \right| \cdot \mu(x) dx=0.$$

 {\bf Definition 2.}  The series (1.2) is said to be universal in weighted spaces $L_\mu^1[0,2\pi]$ in the usual sense, if for any function $f(x) \in L_\mu^1[0,2\pi]$ there exists a growing sequence of natural numbers $n_k$  such that the sequence of partial sums with numbers $n_k$ of the series (1.2) converges to the function $f(x)$ in the metric  $L_\mu^1[0,2\pi]$.

{\bf Definition 3.}  The series (1.2) is said to be universal in
weighted spaces $L_\mu^1[0,2\pi]$ concerning subseries, if for any
function $f(x) \in L_\mu^1[0,2\pi]$ it is possible to choose a
partial series $\displaystyle \sum_{k=1}^\infty f_{n_k}(x)$ from
(1.2), which converges to the $f(x)$ in the metric
$L_\mu^1[0,2\pi]$.

The above mentioned definitions are  given not in the most general
form and only in the generality, in which they will be applied in
the present paper.

In this paper we consider a question on existence of series by
trigonometric system universal in weighted $L_\mu^1[0,2\pi]$
spaces with respect to rearrangements.

Note, that many papers are devoted (see [3]- [10]) to the question
on existence of various types of universal series in the sense of
convergence almost everywhere and on a measure.

Here we will give those results which directly concern to the
Theorems, proved in this paper.

The first usual universal in the sense of convergence almost
everywhere trigonometric series were constructed by D.E.Menshov
[4] and V.Ya.Kozlov [5].

 The series of the form
 $${1\over2}+\sum_{k=1}^\infty a_k cos{kx}+b_k sin{kx} \eqno(1.3)$$
was constructed just by them such that for any measurable on
$[0,2\pi]$  function $f(x)$ there exists the growing sequence of
natural numbers $n_k$  such that the series (1.3) having the
sequence of partial sums with numbers $n_k$   converges to $f(x)$
almost everywhere on $[0,2\pi]$.

Note here, that in this result, when $f(x)\in{L^1_{[0,2\pi]}} $,
it is impossible to replace convergence almost everywhere by
convergence in the metric ${L^1_{[0,2\pi]}}$.

 This result was distributed by A.A.Talalian on arbitrary orthonormal complete systems (see [6]). He also established (see [7]), that if $\{\phi_n(x)\}_{n=1}^\infty $  - the normalized basis of space ${L^p_{[0,1]}},p>1 $, then there exists a series of the form
  $$\sum_{k=1}^\infty{a_k\phi_k(x)},\ \ a_k \to 0.\eqno(1.4)$$
which has property: for any measurable function $f(x)$ the members
of series (1.4) can be rearranged so that the again received
series converge on a measure on [0,1] to $f(x)$.

W. Orlicz [8] observed the fact that there exist functional series
that are universal with respect to rearrangements in the sense of
a.e. convergence in the class of a.e. finite measurable functions.

It is also useful to note that even Riemann proved that every
convergent numerical series which is not absolutely convergent is
universal with respect to rearrangements in the class of all real
numbers.

In [9] it is proved the following result:
\par\par\bigskip

{\bf Theorem 2. } Let $\{ \beta_k\}_{k=0}^\infty$ be a sequence of
positive numbers with $\displaystyle \lim_{k\to \infty}
\beta_k=0$. There exists a trigonometric series of the form
 $$\sum_{k=-\infty}^\infty{C_ke^{ikx} \ \  with \ \ \sum_{k=-\infty}^\infty \left | {C_k} \right|\beta_{|k|} <\infty},\ \ C_{-k}=\overline
{C}_k\eqno(1.5)$$ with the following property: such that for any
number $\epsilon>0$  a weighted function $\mu(x)$, $ 0<\mu(x)
\le1, \left | \{ x\in[0,2\pi]: \mu(x)\not =1 \} \right | <\epsilon
$ can be constructed, so that the series (1.5) is universal in
$L_\mu^1[0,2\pi]$ with respect to rearrangements ( in the usual
sense).
\par\par\bigskip

In this paper we prove the following results.

{\bf Theorem 3.} Let $\omega(t)$ be a continuous function,
increasing in $[0,\infty)$ and $\omega(+0)=0$. Then there exists a
series of the form
 $$\sum_{k=-\infty}^\infty C_ke^{ikx} \ \  with \ \ \sum_{k=-\infty}^\infty C^2_k \omega(|C_k|)<\infty ,\ \ C_{-k}=\overline{C}_k, \eqno(1.6)$$
with the following property: for each $\varepsilon>0$  a weighted
function $\mu(x), 0<\mu(x) \le1, \left | \{ x\in[0,2\pi]:
\mu(x)\not =1 \} \right | <\varepsilon $ can be constructed, so
that the series (1.6) is universal in the weighted space
$L_\mu^1[0,2\pi]$ with respect to rearrangements.

\par\par\bigskip

Analogous of this Theorem for Walsh system was proved by author in
[10].

{\bf Remark.} Using the proofs of Theorem 3 we can construct the
series of the form (1.6) which are universal in the weighted space
$L_\mu^1[0,2\pi]$ with respect simultaneously to rearrangements as
well as to subseries.

The author thanks Professor M.G.Grigorian for his attention to
this paper.

\vskip 4mm

{\bf \S 2. BASIC LEMMA}

\par\par\bigskip
\par\par\bigskip

{\bf Lemma .}  Let $\omega(t)$ be a continuous function,
increasing in $[0,\infty)$ and $\omega(+0)=0$. Then for any given
numbers $0<\varepsilon<{1\over 2}$, $N_0>2$ and a step function
$$f(x)= \sum_{s=1}^q \gamma_s \cdot \chi_{\Delta_s} (x), \eqno(2.1)$$
where  $\Delta_s$ is an interval  of the form $\Delta_m^{(i)}=
\left[ {{i-1}\over {2^m}},{i\over {2^m}} \right] $, $ 1\leq i \leq
2^m $, there exists a measurable set $E \subset [0,2\pi]$ and a
polynomial $P(x)$ of the form
 $$P(x)= \sum_{N_0 \leq |k|<N} C_ke^{ikx} $$
which satisfy the conditions:
$$|E|> 2\pi- \varepsilon, \leqno(1)$$
$$\int_E |P(x)-f(x)|dx<\varepsilon, \leqno(2)$$
$$\sum_{N_0 \leq |k|<N} |C_k|^2\cdot \omega(|C_k|)< \varepsilon,\ \ C_{-k}=\overline {C}_k \leqno(3)$$
$$\max_{N_0 \leq m<N} \left[ \int_e \left | \sum_{N_0 \leq |k|\leq m} C_k e^{ikx} \right | dx \right] <\varepsilon+\int_e |f(x)|dx,\leqno(4)$$
for every measurable subset $e$  of $E$.
\par\par\bigskip
\par\par\bigskip

 {\bf Proof of Lemma .}  Let $0<\epsilon<{1\over 2}$ be an arbitrary number. For any positive number $\eta$ with
$$\eta<{\epsilon^2 \over 4} \cdot \left[ \int_0^{2\pi} f^2(x)dx \right]^{-1},\eqno(2.2)$$
by definition of function $\omega(t)$, there exists a positive
number $\delta<\epsilon$ so that for any $t$, $0<t<\delta$ we have
$$\omega(t)<\omega(\delta)<\eta. \eqno (2.3)$$

 Without restriction of generality, we assume that
$$\max_{1\leq s \leq q}{4 \over \epsilon}|\gamma_s| \cdot \sqrt {|\Delta_s|}<\min \{{\epsilon \over 2},\delta \},\ \ s=1,2,...,q.\eqno(2.4)$$
Set
$$g(x)=\cases {1,\ \ if \ \ x \in [0,2\pi]\setminus \left[ {\varepsilon\cdot \pi \over 2},{3\varepsilon\cdot \pi \over 2}\right] ; \cr   \cr  1-{2 \over\varepsilon},\ \ if \ \ x \in \left[ {\varepsilon\cdot \pi \over 2},{3\varepsilon\cdot \pi \over 2}\right]. \cr} \eqno(2.5)$$
we choose natural numbers $\nu_1$ and $N_1$ so large that the
following inequalities be satisfied:
$${1\over 2\pi}\left|\int_0^{2\pi} g_1(t)e^{-ikt}dt \right|< {\varepsilon \over 16\cdot \sqrt{N_0}},\ \ |k|<N_0, \eqno(2.6)$$
where
$$g_1(x)=\gamma_1\cdot g(\nu_1\cdot x)\cdot \chi_{\Delta_1}(x). \eqno (2.7)$$
(By $\chi_E(x)$ we denote the characteristic function of the set
$E$.)

We put
$$ E_1=\{ x\in \Delta_s:\ \ g_s(x)=\gamma_s \}, \eqno(2.8)$$
By (2.5), (2.7) and (2.8) we have
$$|E_1|>2\pi\cdot (1-\epsilon)\cdot |\Delta_1|;\ \ g_1(x)=0,\ \ x\notin \Delta_1, \eqno(2.9)$$
$$\int_0^{2\pi}g_1^2(x)dx<{2 \over \epsilon}\cdot |\gamma_1|^2\cdot|\Delta_1|. \eqno(2.10)$$

Since the trigonometric system $\{e^{ikx}\}_{k=-\infty}^\infty$ is
complete in $L^2[0,2\pi]$, we can choose a natural number
$N_1>N_0$ so large that

$$\int_0^{2\pi}\left| \sum_{0\leq |k|<N_1}C_k^{(1)}e^{ikx}- g_1(x)\right| dx \leq {\varepsilon \over 8}, \eqno(2.11) $$
where
$$C_k^{(1)}={1\over 2\pi}\int_0^{2\pi} g_1(t)e^{-ikt}dt .$$

Hence by (2.6) and (2.7) we obtain

$$\int_0^{2\pi}\left| \sum_{N_0\leq |k|<N_1}C_k^{(1)}e^{ikx}- g_1(x)\right| dx\leq {\varepsilon \over 8}+\left[ \sum_{0\leq |k|<N_0}
|C_k^{(1)}|^2 \right]^{1 \over 2}< {\varepsilon \over 4}.  $$
\par\par\bigskip

Now assume that the numbers $\nu_1<\nu_2<...\nu_{s-1}$,
$N_1<N_2<...<N_{s-1}$, functions $g_1(x),g_2(x),...,g_{s-1}(x)$
and the sets $E_1,E_2,....,E_{s-1}$ are defined.

We take sufficiently large natural numbers $\nu_s>\nu_{s-1}$ and
$N_s>N_{s-1}$ to satisfy
\par\par\bigskip

$${1\over 2\pi}\left|\int_0^{2\pi} g_s(t)e^{-ikt}dt \right|< {\varepsilon \over 16\cdot \sqrt{N_{s-1}}},\ \ \ \ 1\leq s \leq q,\ \ \ \   |k|<N_{s-1}, \eqno(2.12)$$
$$\int_0^{2\pi}\left| \sum_{0\leq |k|<N_1}C_k^{(s)}e^{ikx}- g_s(x)\right| dx \leq {\varepsilon \over 4^{s-1}}, \eqno(2.13) $$
where
$$g_s(x)=\gamma_s\cdot g(\nu_s\cdot x)\cdot \chi_{\Delta_s}(x),\ \ \ \  C_k^{(s)}={1\over 2\pi}\int_0^{2\pi} g_s(t)e^{-ikt}dt .\eqno(2.14)$$

Set
$$E_s=\{ x\in \Delta_s:\ \ g_s(x)=\gamma_s \}, \eqno(2.15)$$

Using the above arguments (see (2.9)-(2.11)), we conclude that the
function $g_s(x)$ and the set $E_s$ satisfy the conditions:

$$|E_s|>2\pi\cdot (1-\epsilon)\cdot |\Delta_s|;\ \ g_s(x)=0,\ \ x\notin \Delta_s, \eqno(2.16)$$
$$\int_0^{2\pi}g_s^2(x)dx<{2 \over \epsilon}\cdot |\gamma_s|^2\cdot|\Delta_s|. \eqno(2.17)$$
$$\int_0^{2\pi}\left| \sum_{N_{s-1}\leq |k|<N_s}C_k^{(s)}e^{ikx}- g_1(x)\right| dx < {\varepsilon \over 2^{s+1}}.\eqno (2.18)  $$

Thus, by induction we can define natural numbers
$\nu_1<\nu_2<...\nu_q$, $N_1<N_2<...<N_q$, functions
$g_1(x),g_2(x),...,g_q(x)$ and sets $E_1,E_2,....,E_q$ such that
conditions (2.16)- (2.18) are satisfied for all $s,\ \ 1\leq s
\leq q$.
\par\par\bigskip

We define a set $E$ and a polynomial $P(x)$ as follows:

$$E=\bigcup_{s=1}^q E_s, \eqno(2.19)$$

$$P(x)=\sum_{N_0 \leq |k|< N} C_k e^{ikx}=\sum_{s=1}^q \left[ \sum_{N_{s-1}\leq |k|<N_s} C_k^{(s)}e^{ikx} \right], \eqno(2.20)$$

where

$$ C_k=C_k^{(s)} \ \  for\ \ N_{s-1}\leq |k|<N_s,\ \ s=1,2,...,q,\ \ C_{-k}=\overline {C}_k \ \  N=N_q-1. \eqno(2.21)$$

By Bessel's inequality and (2.5), (2.14) for all $s\in [1,q]$ we
get

$$\left[ \sum_{N_{s-1}\leq |k|<N_s} |C_k^{(s)}|^2 \right]^{1 \over 2} \leq \left[ \int_o^{2\pi}g_s^2(x)dx \right]^{1 \over 2}\leq {2 \over \sqrt{\varepsilon}}\cdot |\gamma_s| \cdot \sqrt {|\Delta_s|},\ \ s=1,2,...,q. \eqno (2.22)$$

From (2.5), (2.14) and (2.15) it follows that

$$|E|> 2\pi- \varepsilon. $$
Taking relations (2.1), (2.5), (2.12), (2.14), (2.18) - (2.21) we
obtain
$$\int_E |P(x)-f(x)|dx \leq \sum_{s=1}^q \left[ \int_E \left| \sum_{N_{s-1}\leq |k|<N_s} C_k^{(s)}e^{ikx}-g_s(x) \right| dx\right]<\varepsilon$$
By (2.4), (2.21) and (2.22) for any $k\in [N_0,N]$ we have

$$|C_k|\leq \max_{1 \leq s \leq q}\left[{2 \over  \sqrt {\epsilon}}  \cdot |\gamma_s| \cdot \sqrt {|\Delta_s|}\right]<\delta.$$
From this and (2.3) we get
$$\omega(|C_k|)<\omega(\delta)<\eta,\ \ \forall \ \  k \in [N_0,N].$$
Hence by (2.1), (2.2), (2.4) and (2.22) we obtain
$$\sum_{N_0\leq |k|<N} |C_k|^2\cdot \omega(|C_k|)<\eta \cdot  \sum_{s=1}^q \left[ \sum_{N_{s-1}\leq |k|<N_s} |C_k^{(s)}|^2 \right] < \eta \cdot {4 \over \epsilon} \cdot \left[ \int_0^{2\pi} f^2(x)dx \right]< \epsilon.$$

That is, the statements 1) - 3) of Lemma are satisfied. Now we
will check the fulfillment of statement 4) of Lemma 2.

Let $N_0 \leq m<N$, then for some $s_0,\ \ 1 \leq s_0 \leq q, \ \
\left( N_{s_0} \leq m< N_{s_0+1} \right) $ we will have (see
(2.13) and (2.14))

$$\sum_{N_0 \leq |k|\leq m} C_k e^{ikx}=\sum_{s=1}^{s_0} \left[ \sum_{N_{s-1}\leq |k|<N_s} C_k^{(s)}e^{ikx} \right]+\sum_{N_{s_0-1}\leq |k|\leq m} C_k^{(s_0+1)} e^{ikx}.$$

Hence and from (2.1), (2.4), (2.5), (2.11) and (2.12) for any
measurable set $e \subset E$ we obtain

$$\int_e \left | \sum_{N_0\leq |k|\leq m} C_k e^{ikx} \right | dx \leq$$

$$\leq \sum_{s=1}^{s_0} \left[\int_e \left | \sum_{N_{s-1}\leq |k|<N_s} C_k^{(s)}e^{ikx}-g_s(x) \right|dx \right] +$$

$$+\sum_{s=1}^{s_0} \int_e |g_s(x) | dx +\int_e \left| \sum_{N_{s_0-1}\leq |k|\leq m} C_k^{(s_0+1)} e^{ikx} \right|dx<$$

$$<\sum_{s=1}^{s_0} {\varepsilon \over 2^{s+1}}+\int_e |f(x)|dx+{2 \over \sqrt {\varepsilon}}  \cdot |\gamma_{s_0+1}| \cdot \sqrt {|\Delta_{s_0+1}|} <$$

$$< \int_e |f(x)|dx+ \varepsilon.$$
\par\par\bigskip

{ Lemma is proved.}
\par\par\bigskip

{\bf \S3.PROOF OF THEOREM 3}

\par\par\bigskip
\par\par\bigskip

 Let $\omega(t)$ be a continuous function, increasing in $[0,\infty)$ and $\omega(+0)=0$ and let

$$\{ f_n(x)\}_{n=1}^\infty,\ \  x \in [0,2\pi]\eqno(3.1) $$
\par\par\bigskip
be a sequence of all step functions, values and constancy interval
endpoints of which are rational numbers.Applying Lemma 2
consecutively, we can find a sequence $\{ E_s\}_{s=1}^\infty $ of
sets and a sequence of polynomials

$$P_s(x)=\sum_{N_{s-1}\leq |k|<N_s} C_k^{(s)}e^{ikx} ,\eqno(3.2)$$

$$ 1=N_0<N_1<...<N_s<....,\ \ s=1,2,....,$$

which satisfy the conditions:

$$P_s(x)=f_s(x),\ \ x\in E_s \eqno(3.3)$$

$$\left| E_s\right| >1-2^{-2(s+1)} ,\ \  E_s\subset [0,2\pi],\eqno(3.4)$$

$$ \sum_{N_{s-1}\leq |k|<N_s}\left |C_k^{(s)}\right|\cdot \omega (|C_k^{(s)}|)< 2^{-2s},\ \ C_{-k}^{(s)}=\overline {C}_k^{(s)} \eqno(3.5)$$

$$\max_{N_{s-1}\leq p<{N_s}} \left[ \int_e \left | \sum_{N_{s-1}\leq |k|<p} e^{ikx} \right | dx \right] <2^{-2(s+1)}+\int_e |f_s(x)|dx,\eqno(3.6)$$

for every measurable subset $e$  of $E_s$.

Denote

$$\sum_{k=-\infty}^\infty C_k e^{ikx}=\sum_{s=1}^\infty \left[ \sum_{N_{s-1}\leq |k|<N_s}  C_k^{(s)}e^{ikx} \right],\eqno(3.7)$$

where $ C_k=C_k^{(s)}$ for $N_{s-1}\leq |k|<N_s$, $s=1,2,...$.

Let  $\varepsilon$ be an arbitrary positive number. Setting

$$\Omega_n = \bigcap_{s=n}^\infty E_s,\ \  n=1,2,....;$$

$$ E=\Omega_{n_0} = \bigcap_{s=n_0}^\infty E_s,\ \   n_0=[\log_{1/2} \varepsilon]+1;\eqno(3.8)$$

$$ B=\bigcup _{n=n_0} ^\infty \Omega_n =\Omega_{n_0} \bigcup \left( \bigcup _{n=n_0+1}^ \infty  \Omega_n \setminus \Omega_{n-1} \right). $$

It is clear (see (3.4)) that $\left| B \right|=2\pi$  and $\left|
E \right| >2\pi- \varepsilon .$

We define a function $\mu(x)$ in the following way:

$$\mu(x)=\cases{1\ \ for \ \ x \in E \cup ([0,2\pi] \setminus B);\cr \mu_n \ \ for  \ \ x \in \Omega_n \setminus \Omega_{n-1},\ \ n\geq n_0+1, \cr}\eqno(3.9)$$

where

$$ \mu_n=\left[ 2^{2n}\cdot \prod_{s=1}^n h_s \right]^{-1};\ \ \ \  h_s=|| f_s(x)||_C+ \max_{N_{s-1}\leq p<{N_s}} \Vert  \sum_{N_{s-1}\leq |k|<p}  C_k^{(s)}e^{ikx} \Vert_{C}+1,\eqno(3.10)$$

where

$$ ||g(x)||_C=\max_{x\in [0,2\pi]} |g(x)|,$$

g(x) is a continuous function on  $[0,2\pi]$.

From (3.5),(3.7)-(3.10) we obtain

(A) -- $ \mu(x)$ is a measurable function and
$$ 0<\mu(x) \le1,\ \ \left | \{x\in[0,2\pi]:\mu(x)\not =1\}
\right|<\varepsilon.$$

(B) --   $\displaystyle \sum_{k=1}^\infty |C_k|^2\cdot \omega
(|C_k|)<\infty.$

Hence, obviously we have

 $$\lim_{k\to \infty}{C_k}=0.\eqno(3.11) $$

It follows  from (3.8)-(3.10) that for all $s \geq n_0$  and $p
\in \left[ N_{s-1},N_s \right)$

$$\int_{[0,2\pi] \setminus \Omega_s} \left|  \sum_{N_{s-1}\leq |k|<p}  C_k^{(s)}e^{ikx} \right| \mu(x) dx= $$

$$=\sum_{n=s+1}^ \infty \left[\int_{\Omega_n \setminus \Omega_{n-1}} \left|  \sum_{N_{s-1}\leq |k|<p}  C_k^{(s)}e^{ikx}\right| \mu_n dx \right] \leq $$

$$\leq \sum_{n=s+1}^ \infty2^{-2n} \left[\int_0^ {2\pi} \left|  \sum_{N_{s-1}\leq |k|<p}  C_k^{(s)}e^{ikx}\right| h_s^{-1} dx \right]<{1 \over 3}2^{-2s}.\eqno(3.12) $$

By (3.3), (3.8)-(3.10) for all  $s \geq n_0$ we have

$$\int_0^{1} \left| P_s(x)-f_s(x) \right|\mu(x)dx=\int_{\Omega_s} \left| P_s(x)-f_s(x) \right|\mu(x)dx+$$

 $$+\int_{[0,2\pi] \setminus \Omega_{s}} \left| P_s(x)-f_s(x) \right|\mu(x)dx =2^{-2(s+1)}+$$
$$+\sum_{n=s+1}^\infty \left[\int_{\Omega_n \setminus \Omega_{n-1}} \left| P_s(x)-f_s(x)  \right| \mu_n dx\right] \leq
2^{-2(s+1)}+$$
$$\leq \sum_{n=s+1}^ \infty 2^{-2s}\left[ \int_0^ {2\pi} \left(\left| f_s(x) \right| +\left| \sum_{N_{s-1}\leq |k|<N_s}  C_k^{(s)}e^{ikx} \right| \right) h_s^{-1}dx \right ]<$$
$$<2^{-2(s+1)}+{1 \over 3}2^{-2s}<2^{-2s}.\eqno(3.13) $$

Taking relations (3.6), (3.8)- (3.10) and (3.12) into account we
obtain that for all $p \in \left[ N_{s-1},N_s \right)$ and$s \geq
n_0+1$

$$\int_0^{2\pi} \left| \sum_{N_{s-1}\leq |k|<p}  C_k^{(s)}e^{ikx} \right| \mu(x) dx=$$

$$=\int_{\Omega_s} \left| \sum_{N_{s-1}\leq |k|<p}  C_k^{(s)}e^{ikx} \right| \mu(x) dx+$$

$$+\int_{[0,2\pi] \setminus \Omega_s} \left| \sum_{N_{s-1}\leq |k|<p}  C_k^{(s)}e^{ikx} \right| \mu(x) dx<$$

$$< \sum_{n=n_0+1}^ s \left[\int_{\Omega_n \setminus \Omega_{n-1}} \left| \sum_{N_{s-1}\leq |k|<p}  C_k^{(s)}e^{ikx} \right| dx \right]\cdot \mu_n+{1 \over 3}2^{-2s}<$$

$$<\sum_{n=n_0+1}^ s \left( 2^{-2(s+1)}+\int_{\Omega_n \setminus \Omega_{n-1}} |f_s(x)|dx \right) \mu_n +{1 \over 3}2^{-2s}  =$$

$$=2^{-2(s+1)} \cdot \sum_{n=n_0+1}^ s \mu_n+\int_{\Omega_s} |f_s(x)|\mu(x)dx +{1 \over 3}2^{-2s}<$$

$$<\int_0^{2\pi} |f_s(x)|\mu(x)dx +2^{-2s}. \eqno(3.14)$$
\par\par\bigskip

Let $ f(x) \in L_{\mu}^1 [0,2\pi]$ , i. e.$ \int_0^{2\pi} |f(x)|
\mu(x) dx<\infty$ .
\par\par\bigskip

It is easy to see that we can choose a function $f_{\nu_1}(x)$
from the sequence (3.1) such that

$$\int_0^{2\pi} \left| f(x)- f_{\nu_1}(x) \right|\mu(x)dx<2^{-2},\ \ \nu_1 > n_0+1.\eqno(3.15) $$

Hence, we have

$$\int_0^{2\pi} \left| f_{\nu_1}(x) \right|\mu(x)dx<2^{-2}+\int_0^{2\pi} |f(x)|\mu(x)dx.\eqno(3.16) $$

From (2.1), (A), (3.13) and (3.15)  we obtain with $m_1=1$

$$\int_0^{2\pi} \left| f(x)- \left [ P_{\nu_1}(x)+C_{m_1}e^{im_1x} \right] \right|\mu(x)dx \leq $$

$$\leq \int_0^{2\pi} \left| f(x)- f_{\nu_1}(x) \right|\mu(x)dx+$$

$$+\int_0^{2\pi} \left| f_{\nu_1}(x)-P_{\nu_1}(x) \right|\mu(x)dx+$$

$$+\int_0^{2\pi} \left| C_{m_1}e^{im_1x}\right|\mu(x)dx <2\cdot 2^{-2}+2\pi\cdot\left| C_{m_1}\right|.\eqno(3.17)$$

Assume that numbers
$\nu_1<\nu_2<...<\nu_{q-1};m_1<m_2<...<m_{q-1}$ are chosen in such
a way that the following condition is satisfied:

$$\int_0^{2\pi} \left| f(x)- \sum_{s=1}^j \left [ P_{\nu_s}(x)+C_{m_s}e^{im_sx} \right] \right|\mu(x)dx<$$
$$<2\cdot 2^{-2j}+2\pi \cdot\left|C_{m_j}\right|, \ \ 1\leq j \leq q-1 .\eqno(3.18)$$

We choose a function $f_{\nu_q}(x)$ from the sequence (3.1) such
that

$$\int_0^{2\pi} \left| \left( f(x)- \sum_{s=1}^{q-1} \left [ P_{\nu_s}(x)+C_{m_s}e^{im_sx} \right]  \right)-f_{n_q}(x)\right| \mu(x)dx< 2^{-2q},\eqno(3.19)$$
where $\ \ $$ \nu_q>\nu_{q-1};\ \  \nu_q>m_{q-1}$

This with (3.18) imply

$$\int_0^{2\pi} \left| f_{\nu_q}(x) \right| \mu(x)dx<2^{-2q}+2\cdot 2^{-2(q-1)}+2\pi\cdot\left|C_{m_{q-1}} \right| =$$
$$=9 \cdot 2^{-2q}+ 2\pi\cdot\left|C_{m_{q-1}} \right|.\eqno(3.20) $$
By (3.13), (3.14) and (3.20) we obtain

 $$\int_0^{2\pi} \left| f_{\nu_q}(x)- P_{\nu_q}(x) \right|\mu(x)dx<2^{-2\nu_q},\eqno(3.21)$$
$$\ \ P_{\nu_q}(x)=\sum_{N_{\nu_q-1}\leq |k|<N_{\nu_q}} C_k^{(\nu_q)}e^{ikx}.$$
$$\max_{N_{\nu_q-1} \leq p<N{\nu_q}}  \int_0^{2\pi} \left| \sum_{k=N_{\nu_q-1}}^ p C_k^{(\nu_q)}e^{ikx} \right | \mu(x) dx<10 \cdot 2^{-2q}+2\pi\cdot\left | C_{m_{q-1}} \right |.\eqno(3.22)$$

Denote

$$ m_q= \min \left\{ n \in N: n \notin \left\{  \left\{  \{ k \}_{k=N_{\nu_s-1}}^{N_{\nu_s}-1} \right\}_{s=1}^q \cup \{ m_s\}_{s=1}^{q-1}\right\} \right\}.\eqno(3.23) $$

From (2.1), (A), (3.19) and (3.21) we have
$$\int_0^{2\pi} \left|  f(x)- \sum_{s=1}^q \left [ P_{\nu_s}(x)+C_{m_s}e^{im_sx} \right] \right| \mu(x)dx\leq $$
$$\leq \int_0^{2\pi} \left| \left( f(x)- \sum_{s=1}^{q-1} \left [ P_{\nu_s}(x)+C_{m_s}e^{im_sx} \right]  \right)-f_{\nu_q}(x)\right| \mu(x)dx+$$
$$+\int_0^{2\pi} \left| f_{\nu_q}(x)-P_{\nu_q}(x)\right| \mu(x)dx+$$
$$+\int_0^{2\pi} \left| C_{m_q}e^{im_qx}\right| \mu(x)dx<2 \cdot 2^{-2q}+2\pi\cdot\left| C_{m_q} \right|.\eqno(3.24)$$
Thus, by induction we on $q$ can choose from series (3.7) a
sequence of members
$$C_{m_q}e^{im_qx} ,\ \  q=1,2,...,$$
and a sequence of polynomials
$$P_{\nu_q}(x)=\sum_{N_{\nu_q-1}\leq |k|<N_{\nu_q}}  C_k^{(\nu_q)}e^{ikx},\ \  N_{n_q-1}>N_{n_{q-1}},\ \  q=1,2,....\eqno(3.25)$$

such that conditions (3.22) - (3.24) are satisfied for all $q\geq
1.$

Taking account the choice of $P_{\nu_q}(x)$ and $C_{m_q}e^{im_qx}$
(see (3.23) and (3.25)) we conclude that the series

$$\sum_{q=1}^\infty \left[ \sum_{N_{\nu_q-1}\leq |k|<N_{\nu_q}}  C_k^{(\nu_q)}e^{ikx}+C_{m_q}e^{iqx} \right ]$$

is obtained from the series (3.7) by rearrangement of members.
Denote this series by $\sum C_{\sigma(k)}e^{i\sigma(k) x}.$

It follows from (3.11), (3.22) and (3.24) that the series $\sum
C_{\sigma(k)}e^{i \sigma(k)x}$ converges to the function $f(x)$ in
the metric $L_{\mu}^1[0,2\pi]$, i.e. the series (3.7) is universal
with respect to rearrangements (see Definition 1).

\par\par\bigskip

This completes the proof of Theorem 4.
\par\par\bigskip
\par\par\bigskip

{\bf REFERENCES}

\par\par\bigskip
\par\par\bigskip
\par\par\bigskip

\n [1] N.K.Bary, Trigonometric series, Nauka, Moscow,1961; English
trans. in Pergamon Press, Oxford, 1964.
\par\par\bigskip

\n [2] K.S.Kazarian, R.Zink ,Some ramifications of a theorem of
Boas and Pollard concerning the completion of a set of functions
in $L^2$, Trans.Amer.Math.Soc., v.349, n.11, p. 4367-4383.
\par\par\bigskip

\n [3] M. G. Grigorian "On the representation of functions by
orthogonal series in weighted $L^p$ spaces,\ \ Studia. Math.
134(3)1999, 211-237.
\par\par\bigskip

\n [4] D.E.Menshov, On the partial summs of trigonometric series,
Mat. Sb. 20(1947), 197-238 [in russian].
\par\par\bigskip

\n [5] V.Ya.Kozlov, On the complete systems of orthogonal
functions, Mat. Sb. 26(1950), 351-364 [in russian].
\par\par\bigskip

\n [6] A.A.Talalian, On the convergence almost everywhere the
sumbsequence of partial summs of general orthogonal series,Izv.
Ak. Nauk Arm. SSR ser. Math. 10(1957), 17-34 [in russian].
\par\par\bigskip

\n [7] A.A.Talalian, On the universal series with respect to
rearrangements, Izv. AN. SSSR ser. Math. 24(1960), 567-604 [in
russian].
\par\par\bigskip

\n [8] W.Orlicz , Uber die unabhangig von der Anordnung fast
uberall kniwergenten Reihen, Bull. de l'Academie Polonaise des
Sciences, 81 (1927), p. 117-125.
\par\par\bigskip

\n [9]  S.A.Episkoposian, M.G.Grigorian,"On universal
trigonometric series in weighted spaces $L_\mu^1[0,2\pi]$ " , East
Journal on Approximations, 1999,v.5 , n.4, 483-492.
\par\par\bigskip

\n [10] S.A. Episkoposian , "On the existence of universal series
by Walsh system ",  Izvestiya Natsionalnoi Akademii Nauk Armenii,
English trans. in: Journal of Contemporary Mathematical Analysis,
2003, v. 38  , n.4, p.25-40.

\par\par\bigskip
\par\par\bigskip

Department of Physics,

State University of Yerevan,

Alex Manukian 1, 375049 Yerevan, Armenia

e-mail: sergoep@ysu.am

\end{document}